%% file: lincei11900.tex
\documentstyle[epsfig,eqsection]{article}

\include{macro}

\begin{document}

\title{{\bf Diffusion time and splitting of separatrices for
nearly integrable isochronous Hamiltonian systems}}

\author{Massimiliano Berti and Philippe Bolle}
\date{}
\maketitle

{\bf Abstract:}
We consider the problem of Arnold's diffusion for
nearly integrable isochronous Hamiltonian systems.
We prove a shadowing theorem which improves the known
estimates for the diffusion time.
We also justify for three time scales systems
that the splitting of the separatrices
is correctly predicted by the Poincar\'e-Melnikov function.
\footnote{Supported by M.U.R.S.T. Variational Methods and Nonlinear
Differential Equations.}
\\[2mm]
\\[2mm]
\indent
{\bf Riassunto:} {\it Tempo di diffusione e splitting delle separatrici
per sistemi Hamiltoniani isocroni quasi-integrabili}.
Consideriamo il problema della diffusione di Arnold per
sistemi Hamiltoniani isocroni quasi-integrabili. Dimostriamo un teorema di
shadowing che migliora le stime sul tempo di diffusione sinora note.
Giustifichiamo inoltre, per sistemi a tre scale temporali, che
lo splitting delle separatrici \`e correttamente previsto
dalla funzione di Poincar\'e-Melnikov.



\section{Introduction}

We outline in this Note some recent results on Arnold's diffusion for
nearly integrable isochronous Hamiltonian systems: complete
proofs and further results are contained in \cite{BB1}.
We consider Hamiltonians $ {\cal H}_\mu $ of the form
\be\label{eq:Hamg}
{\cal H}_\mu = \om \cdot I + \frac{p^2}{2} + (\cos q - 1) +
\mu f( \vphi, q )
\ee
with angles variables  $ ( \vphi, q ) \in {\bf T}^n \times {\bf T}^1 $
and action variables $ ( I, p ) \in {\bf R}^n \times {\bf R}^1$. We assume
\begin{itemize}
\item
$(H1)$ There exists $ \gamma > 0 $, $ \tau > n $
such that $| \om \cdot k| \geq \gamma / |k|^{\tau} $,
$ \forall k \in {\bf Z}^n, k \neq 0$.
\end{itemize}
\noindent
Hamiltonian $ {\cal H}_\mu $ describes a
system of $n$ isochronous harmonic oscillators
of frequencies $ \om $ weakly coupled with a pendulum.
When $ \mu = 0 $ the energy
$ \om_i I_i $ of each oscillator  is a
constant of the  motion.
The problem of {\it Arnold's diffusion} is whether,
for $ \mu \neq 0 $, there exists motions whose net effect is to transfer
energy from one oscillator to the others.

Stemming from \cite{Arn} these kind of results are usually proved noting
that,
for $ \mu =0 $ Hamiltonian ${\cal H}_\mu$ admits
a continuous family of
$ n $-dimensional partially hyperbolic invariant tori
${\cal T}_{I_0} = \{ (\vphi , I, q,p) \in  {\bf T}^n
\times  {\bf R}^n \times {\bf T}^1 \times {\bf R}^1 \ |
\  I = I_0, \ q = p = 0 \}$
possessing stable and unstable manifolds
$ W^s ({\cal T}_{I_0}) = W^u({\cal T}_{I_0}) =
\{ (\vphi , I, q,p) \in  {\bf T}^n
\times  {\bf R}^n \times {\bf T}^1 \times {\bf R}^1 \  |
  \ I = I_0, \ p^2 / 2 + ( \cos q - 1) = 0 \}.$
By $(H1)$ all the unperturbed tori ${\cal T}_{I_0}$, with their
stable and unstable manifolds, persist, for $ \mu  $ small enough,
being just sligthly deformed. The perturbed stable and unstable manifolds
$W^s_\mu ({\cal T}_{I_0}^\mu )$ and
$W^u_\mu ({\cal T}_{I_0}^\mu) $ may split and intersect transversally
giving rise to a chain of tori connected by heteroclinic orbits.
By a shadowing type argument one can then prove the existence of an orbit
such that the action variables $ I $ undergo a variation of $ O (1) $
in a certain time $ T_d $ called the {\it diffusion time}.
In order to prove the existence of diffusion orbits
following the previous mechanism one encounters two different problems

\begin{description}
\item{$1)$} Shadowing theorem;

\item{$2)$} Splitting of separatrices;
\end{description}

By means of a variational technique inspired from
\cite{BB}-\cite{BBS}-\cite{AB}
we give in the next section a general shadowing theorem which
improves, for isochronous systems, the known estimates
on the diffusion time obtained
in \cite{CG}-\cite{G1}-\cite{M}-\cite{C1} by geometric methods
and in \cite{BCV} by Mather's theory.
In section 3, using methods introduced in \cite{An1}, we give some
results on the splitting of the separatrices.

\section{The shadowing  theorem}

\subsection{ Perturbation term vanishing
along the invariant tori}

We first describe our results when the perturbation term is
$ f( \vphi, q) = (1 - \cos q ) \ f(\vphi )$
so that the tori $ {\cal T}_{I_0} $ are still invariant
for $\mu \neq 0$.
The equations of motion derived by Hamiltonian ${\cal H}_\mu$ are
\be\label{eqmotion}
\dot{\vphi} = \om, \qquad
\dot{I} = - \mu (1 - \cos q ) \ \partial_{\vphi}f(\vphi ), \qquad
\dot{q} = p, \qquad
\dot{p} = \sin{q} - \mu \sin q \ f (\vphi ).
\ee
The angles $ \vphi $ evolve  as
$ \vphi (t) = \om t + A $ and then (\ref{eqmotion})
can be reduced  to the quasi-periodically forced pendulum
equation
\be\label{pendper}
- \ddot{q} + \sin{q} = \mu \ \sin q \ f ( \om t + A),
\ee
corresponding to the Lagrangian
\be\label{lagraper}
{\cal L}_\mu (q, \dot{q},t) = \frac{{\dot q}^2}{2} + (1- \cos q) +
\mu (\cos q - 1)f(\om t + A).
\ee
For each solution $ q(t) $ of (\ref{pendper}) one recovers the dynamics
of the actions $ I (t) $ by quadratures in (\ref{eqmotion}).
\\[1mm]
\indent
For $ \mu = 0 $ equation (\ref{pendper})
possesses the family of heteroclinic solutions
$ q_\teta (t) = 4 \ {\rm arctg} (\exp {(t - \teta)}),
\  \teta \in {\bf R} $.
Using the Contraction Mapping Theorem as in \cite{BB} one can prove that
near the unperturbed heteroclinic solutions $ q_\teta (t) $
there exist, for $ \mu $ small enough, {\it
``pseudo-heteroclinic solutions''}
$ q_{A,\teta}^\mu (t)$ of equation (\ref{pendper}).
$ q_{A,\teta}^\mu (t)$ are true solutions of (\ref{pendper})
in each interval $(- \infty, \teta)$ and $(\teta, +\infty)$;
at time $t = \teta $ such pseudo-solutions are glued with continuity
at value  $ q^\mu_{A, \teta} ( \teta ) = \pi $ and for
$ t \to \pm \infty $ are asymptotic
to the equilibrium $0 $ mod $ 2 \pi $.
Moreover, by a uniqueness property,
$q^{\mu}_{A,\theta}$ depends smoothly on
$(\mu , A , \theta)$.
We can then define the function
$ F_\mu : {\bf T}^n \times {\bf R} \to {\bf R} $
as the action functional of Lagrangian (\ref{lagraper})
evaluated on the ``1-bump pseudo-heteroclinic solutions''
$q_{A, \teta }^\mu (t)$, namely
\be\label{def:homo}
F_\mu (A, \teta) =
\int_{ - \infty}^\teta {\cal L}_\mu ({  q}_{A,\teta}^\mu (t),
\dot{  q}^\mu_{A,\teta} (t), t) \ dt
+ \int_\teta^{+ \infty}
{\cal L}_\mu ({  q}_{A,\teta}^\mu (t),
\dot{  q}^\mu_{A, \teta } (t),t) \ dt,
\ee
and the {\it ``homoclinic function''} $ G_\mu : {\bf T}^n \to {\bf R} $ as
\be\label{eq:homo}
G_\mu ( A ) = F_\mu ( A, 0 ).
\ee
There holds $ F_\mu  (A , \teta ) = G_\mu (A + \om \teta), \forall
\teta \in {\bf R}.$

\begin{remark}
The homoclinic function $ G_\mu $ is the difference
between the generating functions
$ {\cal S}_{\mu,I_0}^\pm (A, q )$ of the stable and the unstable manifolds
$ W_\mu^{s,u} ( {\cal T}_{I_0} )$
(which in this case are {\it exact} Lagrangian manifolds)
at the  section $ q = \pi $, namely
$G_\mu (A) = {\cal S}_{\mu, I_0}^- (A, \pi ) -
{\cal S}_{\mu, I_0}^+ ( A, \pi ) $. Note that
$G_\mu$ is independent of $I_0$.
\end{remark}

We now give an example
of  condition on $G_{\mu}$ which implies the existence
of diffusion orbits.

\begin{condition}\label{spli}
({\it ``Splitting condition''})
There exist $A_0 \in {\bf T}^n $, $ \delta > 0$, $ 0 < \a < \rho $ such that
\begin{itemize}
\item $(i)$
$\inf_{\partial B_\rho (A_0)} G_\mu > \inf_{B_\rho ( A_0 )} G_\mu + \delta$;

\item $(ii)$
$ \sup_{ B_\alpha (A_0) } G_\mu < \frac{\delta}{4}
+ \inf_{B_\rho (A_0)} G_\mu $;
\item $(iii)$
$d(\{ A \in B_\rho ( A_0 ) \ | \
G_\mu (A)  < \delta / 2  + \inf_{B_\rho ( A_0 )} G_\mu  \},
\{ A \in B_\rho ( A_0 ) \ | \ G_\mu (A) \geq 3\delta / 4
+  \inf_{B_\rho ( A_0 )} G_\mu \})
\geq 2 \alpha $.
\end{itemize}
\end{condition}

Note that the above ``splitting condition''
is clearly satisfied if $ G_\mu $ possesses in
$ A_0 \in {\bf T}^n $ a non-degenerate minimum.
Moreover $B_\rho (A_0)$,  open ball of radius $\rho$,
could be replaced by an open subset $U$ of ${\bf T}^n$ whose covering
set is bounded in ${\bf R}^n$.

The following shadowing type theorem holds

\begin{theorem}\label{thm:main}
Assume $ ( H1 ) $ and let $ G_\mu $ satisfy
the splitting condition \ref{spli}.
Then $ \forall I_0, I_0'$ with $ \om \cdot I_0 = \om \cdot I_0'$,
there is a heteroclinic orbit connecting the invariant tori
${\cal T}_{I_0}$ and ${\cal T}_{I_0'}$.
Moreover, there exists $ C > 0 $ such that
$ \forall  \eta > 0 $ small  the ``diffusion time'' $T_d$ needed
to go from a $\eta$-neighbourhood of
${\cal T}_{I_0}$ to a $\eta$-neighbourhood of
${\cal T}_{I_0'}$ is bounded by
\be\label{timediff}
T_d \leq C \frac{| I_0 - I_0' |}{\delta} \rho \max \Big( | \ln \delta |,
\frac{1}{\alpha^\tau} \Big) + C | \ln ( \eta )|.
\ee
\end{theorem}

\begin{remark}
The meaning of (\ref{timediff}) is the following:
the diffusion time $ T_d $ is estimated by the product
of the number of heteroclinic transitions
$ k = $ ( heteroclinic jump / splitting ) $ = | I_0' - I_0 | / \delta $,
and of  the time $ T_s $ required for a single transition,
that is $ T_d = k \cdot T_s  $.
The time for a single transition $ T_s $ is bounded by the maximum
time between the
``ergoditation time'' $ ( 1 / \alpha^\tau ) $, i.e.
the time needed for the flow $ \om t $ to make an
$ \alpha $-net of the torus,
and the time $ | \ln \delta | $ needed  to ``shadow''
homoclinic orbits for the forced pendulum equation. We use
here that these homoclinic orbits are exponentially
asymptotic to the equilibrium.
\end{remark}

When the frequency vector $ \om $ is considered as a
constant, independent of any parameter (``a-priori unstable'' case),
it is easy to give a criterion for
the splitting condition \ref{spli} thanks to the first order 
expansion in $ \mu $. There holds
\be
G_\mu (A )  = const + \mu \Gamma (A) + O( \mu^2 ),
\ee
where $ \Gamma: {\bf T}^n \to {\bf R} $ is the Poincar\'e-Melnikov primitive
\be
\Gamma (A) =
\int_{\bf R} ( 1 - \cos q_0 (t) ) f (\om t + A) \ dt.
\ee
As a corollary of theorem \ref{thm:main} we have

\begin{corollary}\label{thm:apu}
Assume $ (H1) $ and let $ \Gamma $ possess a
non-degenerate minimum. Then, for $ \mu $ small enough,
the  statement of theorem \ref{thm:main} holds where
the diffusion time is
\be\label{eq:estiapu}
T_d = O (\frac{1}{\mu} \log { \frac{1}{\mu } } ).
\ee
\end{corollary}

\begin{remark}
The estimate on the diffusion time obtained
in \cite{CG} is $ T_d >> O( \exp{1/\mu} )$
and is improved in \cite{G1} to be $ T_d = O(\exp{1/ \mu}) $.
Recently in \cite{BCV} by means of Mather's theory
the estimate on the diffusion time has been improved to be
$T_d = O ( 1/ \mu^{2 \t +1} ) $.
In \cite{C1} it is obtained via geometric methods
that $ T_d = O ( 1/ \mu^{\t +1} ) $.
The main reason for which we manage to improve also the estimates
of  \cite{BCV} and \cite{C1} is that
the shadowing orbit of theorem \ref{thm:main}
can be chosen, at each transition, to approach the
homoclinic point $ A_0 $, only up to the distance $ \a $ which does
not depend of $ \mu $.
Note moreover that estimate (\ref{eq:estiapu}) is independent
of the number of rotators $ n $.
\end{remark}

\begin{remark}
The above result answers to a question raised in
\cite{Lo} (sec.7) proving that, at least for isochronous systems,
it is possible to reach the maximal speed of diffusion
$\mu / | \log \mu |$ (moreover independently on the dimension $ n $).
\end{remark}

\subsection{More general perturbation term}

Dealing with more general perturbations
$f(\vphi, q)$ the first step is to prove
the persistence of invariant tori for $ \mu \neq 0 $ small enough.
It appears that no more than the standard
Implicit Function Theorem
is required to prove the following
well known result.

\begin{theorem}\label{thm:invto}
Let $ \om $ satisfy $(H1)$. For $ \mu $ small enough,
for all $ I_0 \in {\bf R}^n $
system ${\cal H}_\mu $ possesses $ n $-dimensional invariant tori
$ {\cal T}_{I_0}^\mu \approx {\cal T}_{I_0} $ of the form
$$
{\cal T}_{I_0}^\mu = \Big\{
I = I_0 + a^\mu ( \psi ), \ \vphi = \psi, \
q =  Q^\mu ( \psi ), \
p =  P^\mu ( \psi ), \quad \psi \in {\bf T}^n \Big\},
$$
with $ Q^\mu (\cdot ), P^\mu (\cdot), a^\mu ( \cdot ) =  O(\mu )
$.
The dynamics on ${\cal T}_{I_0}^\mu$ is conjugated
to the rotation of speed $ \om $.
\end{theorem}

In order to reduce to the previous case we want to put the tori ${\cal
T}^\mu_{I_0}$ ``at the origin'' in the $(q,p)$ coordinates
by a symplectic change of variables.
As ${\cal T}_{I_0}^\mu $ is  isotropic, the transformation of
coordinates
$(\psi,J, u,v) \to ( \vphi, I,  q, p)$
defined, on the covering space $ {\bf R}^{2 ( n + 1 )} $ of
$ {\bf T}^n  \times {\bf R}^n \times {\bf T}^1 \times
{\bf  R}^1 $, by
$$
\vphi = \psi, \quad
I =  a^\mu (\psi ) + \partial_{\psi} P^\mu (\psi) \cdot u -
\partial_{\psi} Q^\mu (\psi) \cdot v + J, \quad   q  = Q^\mu (\psi) + u,
\quad  p  = P^\mu (\psi) + v
$$
is symplectic.
In the new coordinates each invariant torus ${\cal T}_{I_0}^\mu $
is simply described by
$  \ \psi \in {\bf T}^n, \ J=I_0,  \  u = v = 0 $ and
the new Hamiltonian writes
$$
{\cal K}_\mu = E_\mu + \om \cdot J + \frac{v^2}{2} + ( \cos u - 1)  +
P_0 (\mu, u, \psi) \leqno{({\cal K}_{\mu})}
$$
where the perturbation term is
$$
P_0 (\mu,  u, \psi ) = \Big( \cos( Q^\mu + u) - \cos Q^\mu +
(\sin Q^\mu ) u + 1 - \cos u \Big) +
\mu \Big( f( \psi, Q^\mu  + u) -  f ( \psi, Q^\mu ) -
\partial_q f ( \psi, Q^\mu   ) u \Big)
$$
and  $ E_\mu $ is the energy of the perturbed
invariant torus $ ( a^\mu ( \psi ), \psi, Q^\mu (\psi ), P^\mu (\psi )).$
Hamiltonian ${\cal K}_\mu$  corresponds to
the quasi-periodically forced pendulum
\be\label{pendtrasf}
- \ddot{u} + \sin{u} =  \partial_u P_0 (\mu, u,  \om t + A)
\ee
of Lagrangian
\be\label{eq:lagrap}
L_\mu = \frac{\dot{u}^2}{2} + ( 1- \cos u) - P_0 (\mu, u, \om t + A).
\ee
Since the Hamiltonian $ {\cal K}_\mu $ is not periodic in the variable $ u $
we cannot directly apply theorem \ref{thm:main} and
the arguments of the previous section require some modifications.
For $ \mu $ small enough there exists, near $ q_\teta $
(or more exactly its covering orbit in ${\bf R}$),
a unique {\it pseudo-heteroclinic} solution
$ u_{A, \teta}^\mu (t) $, which satisfies (\ref{pendtrasf})
in $(- \infty, \teta)$ and $(\teta, +\infty)$, is
glued with continuity at value  $ u^\mu_{A, \teta} ( \teta ) = \pi $
and is asymptotic as $ t \to -\infty $ (resp. $+\infty$) to the equilibrium
$0
$ (resp. $2\pi$).
Then we define the function ${\cal F}_\mu ( A, \theta) $ as
\begin{eqnarray*}
{\cal F}_\mu ( A, \theta) & = & \int_{-\infty}^\theta
\frac{({\dot u}^\mu_{A, \theta})^2}{2} + ( 1 - \cos u^\mu_{A, \theta} )
- P_0 (\mu, u^\mu_{A, \theta} , \omega t + A) \ dt \\
& + & \int_\theta^{+\infty}
\frac{({\dot u}^\mu_{A, \theta})^2}{2} + ( 1 - \cos u^\mu_{A, \theta} )
- P_1 ( \mu, u^\mu_{A, \theta}, \om t + A) \ dt
+ 2 \pi {\dot q}^\mu_A ( \theta ),
\end{eqnarray*}
where $ q^\mu_A (t) = Q^\mu (\om t + A) $ and
\begin{eqnarray*}
P_1 (\mu, u, \omega t + A) & = & \Big( \cos (q^\mu_A (t)+ u ) -
\cos q^\mu_A (t)  + \sin q^\mu_A (t) \ ( u - 2 \pi ) + 1 - \cos u \Big) \\
& + &
\mu \Big( f(\omega t + A, q^\mu_A (t) + u) - f( \omega t + A, q^\mu_A
(t) ) -
(\partial_q f) (\omega t + A, q^\mu_A (t) ) \ (u - 2 \pi ) \Big).
\end{eqnarray*}
We define the {\it ``homoclinic function''} $ {\cal G}_\mu ( A )$ as
${\cal G}_\mu ( A ) = {\cal F}_\mu ( A, 0 )$.
\\[1mm]
The term $ 2 \pi {\dot q}^\mu_A ( \theta ) $  takes into account
the fact that the stable and the unstable manifolds
of the tori $ {\cal T}^\mu_{I_0} $ are
{\it not exact} Lagrangian manifolds, see \cite{LMS}.
We have

\begin{theorem}\label{thm:main1g}
Assume $(H1)$ and let ${\cal G}_\mu$ satisfy condition \ref{spli}.
Then the statement of theorem \ref{thm:main} holds.
\end{theorem}
For $ \mu $ small enough we have
\be\label{eq:pm}
{\cal G}_\mu (A ) = const + \mu M(A) + O(\mu^2), \qquad \forall A \in
{\bf T}^n
\ee
where
$ M ( A ) = \int_{-\infty}^{+\infty}
\Big[ f ( \omega t + A, q_0 (t)) - f ( \omega t + A, 0 ) \Big] \ dt$.

\begin{corollary}\label{apug}
Assume $(H1)$ and let $ M $ possess a non-degenerate minimum.
Then, for $\mu$ small enough, there exists a diffusion orbit
with diffusion time $T_d = O (1/ \mu \log (1/ \mu ) )$.
\end{corollary}

\section{Splitting of separatrices}

\subsection{Approximation of the homoclinic function}

If the frequency vector $ \om = \om_\e $ contains some
``fast frequencies'' $ \beta_i / \e^b $, $ b > 0 $,  $\e$
being a small parameter,  the oscillations of
the Melnikov function along some directions turn out to be exponentially
small with respect to $ \e $.
Hence the development  (\ref{eq:pm}) will
provide a valid measure of the splitting only for
$ \mu $ exponentially small with respect to $ \e $.
In order to justify the dominance
of the Poincar\'e-Melnikov function when $ \mu = O(\e^p )$
we need more refined estimates for the error.
The classical way to overcome this difficulty would be to extend
analytically the function $ F_\mu (A, \teta) $ for
complex values of the variables, see \cite{LMS}-\cite{DGTJ}.
However it turns out that the function $F_\mu (A, \teta) $ can {\it
not} be easily analytically extended in a 
sufficiently wide complex strip (roughly speaking, the condition
$ q^\mu_{A, \teta } ( Re \ \teta ) = \pi $ appearing
naturally when we try to extend the definition of 
$q^\mu_{A, \teta}$ to $\theta \in {\bf C}$ breaks
analyticity). We bypass this problem considering the action functional
evaluated on different ``1-bump pseudo-heteroclinic solutions''
$ Q^\mu_{A,\teta} $.
This new ``reduced action functional'' $ {\wtilde F}_\mu (A, \teta) $
has the advantage to have an analytical extension defined
for $\theta \in {\bf R}+ i (-\pi/2 , \pi/2 )$.
More precisely let us assume that
$ f( \vphi, q ) = ( 1 - \cos q ) f ( \vphi )$ and that
$ f $ can be extended to an analytical function over
$D:=  ( {\bf R} + i [ - a_1, a_1 ]) \times \ldots \times
({\bf R} + i [ -a_n, a_n ])$, for some $ a_i \geq 0 $.
Then
$$
f( \vphi  ) = \sum_{k \in {\bf Z}^n} f_k \exp^{i k \cdot \vphi}
\quad {\rm with} \quad
| f_k | \leq \frac{ C_s}{|k|^s} \exp \Big( - \sum_{i=1}^n a_i |k_i| \Big),
$$
for all $s \in {\bf N}$.

Define $\psi_0: {\bf R} \to {\bf R} $ by
$ \psi_0 (t) = \cosh^2 ( t )/( 1 + \cosh t)^3 $ and
set $\psi_\teta (t)= \psi_0 (t - \teta)$.
Note that $ \int_{\bf R} \psi_0 (t) {\dot q}_0 (t) \ dt \neq 0 $.
By the Contraction Mapping Theorem we find
near $ q_\teta $, for $ \mu $ small enough,
{\it pseudo-heteroclinic solutions} $ Q_{A, \teta}^\mu (t) $
and a constant $ \a_{A, \teta}^\mu $ defined by
$$
- \ddot{Q}_{A, \teta}^\mu + \sin{Q}_{A, \teta}^\mu = \mu \
\sin Q_{A, \teta}^\mu f (\om t + A) +\a_{A, \teta}^\mu \psi_\teta (t) \quad
{\rm and} \quad
\int_{\bf R} \Big( Q_{A, \teta}^\mu (t)- q_\teta (t) \Big)
\psi_\teta (t) \ dt = 0.
$$
We define the function
$\wtilde{F}_\mu: {\bf T}^n \times {\bf R} \to {\bf R} $
as the action functional of Lagrangian  (\ref{lagraper})
evaluated on the ``1-bump pseudo-heteroclinic solutions''
$ Q_{A, \teta }^\mu (t) $, namely
$$
\wtilde{F}_\mu (A, \teta) =
\int_{\bf R}  {\cal L}_\mu ( Q_{A,\teta}^\mu (t),
\dot Q^\mu_{A,\teta} (t), t) \ dt
$$
and $ \wtilde{G}_\mu: {\bf T}^n  \to {\bf R} $ as
$ \wtilde{G}_\mu ( A ) = \wtilde{F}_\mu ( A, 0 ).$
\noindent
The relation between the functions $ \wtilde{G}_\mu $ and $ G_\mu $ is
given below

\begin{theorem}\label{thm:diffeo}
There exists a smooth diffeomorphism
$ \psi_\mu : {\bf T}^n \to {\bf T}^n $ of the form
$ \psi_\mu (A) = A + g_\mu (A) \om $ with
$ g_\mu ( A ) : {\bf T}^n \to {\bf R} $ satisfying
$ (g_\mu ( A ), \partial_A g_\mu ( A ) ) \to O $ as $ \mu \to 0 $,
such that
$ G_\mu  = \wtilde{G}_\mu \circ \psi_\mu. $
\end{theorem}

We now approximate the Fourier coefficients
of $ \wtilde{G}_\mu (A) = \sum_{k \in {\bf Z}^n }
\wtilde{G}_k \exp^{i k \cdot A  }$ with the
Fourier coefficients of the Poincar\'e-Melnikov primitive
$ \Gamma (A) = \sum_{k \in {\bf Z}^n } \Gamma_k \exp^{i k \cdot A}$.
$ \Gamma_k $ are explicitely given by
\be\label{eq:melcoe}
\Gamma_k = f_k \frac{2 \pi (k \cdot \om)}{\sinh (k \cdot \om
\frac{\pi}{2})} \ee
Set $ || f || = \sup_{A \in  D} | f ( A ) |.$

Note that the unperturbed separatrix
$ q_0 (t) = 4 \ {\rm arctg} (\exp ( t ))$ can be
analytically extended up to $ |Im \ t | < \pi / 2 $.
Using  the Contraction Mapping Theorem,
it is possible to extend analytically the
function $ \wtilde{F}_\mu ( A , \teta ) $
up to the strip $ D \times ({\bf R} +
i(- \pi / 2 + \d, \pi / 2 - \d))$,
provided $ \mu || f || \delta^{-3} $ is small.
By an estimate of   $\wtilde{F}_\mu ( A , \teta )-\mu \Gamma
(A+\om \theta) $ over its complex domain
and a standard lemma on Fourier coefficients of analytical functions we
obtain

\begin{theorem}
There is a constant $C$ such that, for
$ \mu ||f|| \d^{-3} $  small enough,
$ \forall k \neq 0, k \in {\bf Z}^n $,
$\forall \d \in (0, \frac{\pi}{2}) $,
\be\label{eq:devecoe}
| \wtilde{G}_k - \mu \Gamma_k | \leq
\frac{C \mu^2 || f ||^2 }{ \d^4}
\exp \Big( - \sum_{i=1}^n a_i | k_i | \Big)
\exp \Big( - | k \cdot \om | ( \frac{\pi}{2} - \d) \Big).
\ee
\end{theorem}

\subsection{Three time scales}

We consider three time scales Hamiltonians of the form
$$
{\cal H} = \frac{I_1}{ \sqrt{ \e } } + \e^a I_2 +\frac{ p^2 }{2}
+ ( \cos q - 1 ) + \mu ( \cos q - 1 ) f( \vphi_1, \vphi_2), \
I_1 \in {\bf R}, I_2 \in {\bf R}^{n-1}, \ n \geq 2,
$$
namely ${\cal H}_\mu $ with
$ \om_\e = ( \frac{1}{\sqrt{ \e }}, \e^a ) $.
Such systems have been dealt with for example in
\cite{GGM} and \cite{PV}.

Let $ \mu \e^{- 3 / 2 } $ be small enough. We
assume only that  $ f $ is analytical w.r.t $ \vphi_2 $ (more
precisely, $a_1=0$, and for $i\geq 2$, $a_i>0$, $a_i > \pi /2$ if $a=0$).
Set
$ \wtilde{G}_\mu (A) = \sum_{k_1 \in {\bf Z} }
\wtilde{G}_{k_1} ( A_2 ) \exp^{i k_1 \cdot A_1 }$ and
$ \Gamma (A) = \sum_{k_1 \in {\bf Z}}
\Gamma_{k_1} (A_2) \exp^{i k_1 \cdot A_1  }$.
From estimate (\ref{eq:devecoe}) we obtain
\begin{theorem}\label{thm:tts}
For $\mu \e^{- 3/2} $ small there holds
\begin{eqnarray*}
\wtilde{G}_\mu (A_1, A_2) & = & Const +
\Big( \mu \Gamma_0 (\e, \mu, A _2) + R_0 (\e, \mu, A_2) \Big) +
2 {\rm Re} \
\Big( \mu \Gamma_1 (\e, \mu, A _2)
+ R_1 (\e, \mu, A_2) \Big) e^{i A_1}\\
& + &
O (\mu  \ep^{-1/2} ||f|| \exp^{ - \frac{\pi}{\sqrt{\e}} }  )
\end{eqnarray*}
where
$$
R_0 (\e, \mu,  A_2) = O \Big( \mu^2 ||f||^2 \Big) \quad {\rm and}
\quad R_1 (\e, \mu, A_2) = O \Big( \frac{\mu^2 ||f||^2}{\e^2}
\exp^{- \frac{\pi}{2 \sqrt{\e}}} \Big).
$$
\end{theorem}

\begin{remark}
(i) This improves the results in \cite{PV} which require
$ \mu = \e^p $ with $ p > 2 + a $.

(ii) Note that theorem \ref{thm:tts} certainly holds in any dimension,
while the results of \cite{GGM}, which hold for more general systems,
are proved for 2 rotators only.

(iii) Theorem \ref{thm:tts} is not in contradiction
with \cite{GGM3}.
\end{remark}

This theorem jointly with theorem \ref{thm:diffeo} enables us to provide
conditions implying the existence of diffusion orbits. In fact,
if $\wtilde{G}_{\mu}$ has a proper minimum satisfying
condition \ref{spli}, so has ${G}_{\mu}$.
For instance we obtain the following result

\begin{theorem} \label{lem:3ts}
Assume that there are $\ov{A}_2$ and $d,c >0$ such that,
for all small $\e >0$,
$|\Gamma_1 ( A_2 )| > (c/\sqrt{\e}) e^{-\pi /( 2 \sqrt{\e})} $
for all $|A_2 - \overline{A}_2| < d$ and
$\Gamma_0 (\ov{A}_2 \pm d) > \Gamma_0 (\ov{A}_2) +c$.
Then, for $\mu \e^{-3/2}$ small enough, condition \ref{spli}
is satisfied by ${G}_{\mu}$, with
$ \alpha = \ov{ C} e^{-\pi / (2 \sqrt{\e})}$, $\delta =
c\mu/(2\sqrt{\e})
e^{-\pi / (2 \sqrt{\e})}$, $\ov{C}$ constant.
\end{theorem}

\begin{remark}
In order to prove the
splitting of the separatrices using theorem \ref{thm:tts}
it is necessary, according with \cite{GGM} and \cite{PV},
that $ \exists m, l \in {\bf Z}^{n-1}$ such that
$ f_{0,l}, f_{1, m} \neq 0 $.
\end{remark}

\textbf{Acknowledgments}
The first author wishes to thank Prof. G. Gallavotti
for stimulating discussions.

\noindent
{\it Massimiliano Berti, S.I.S.S.A., Via Beirut 2-4,
34014, Trieste, Italy, berti@sissa.it}.
\\[2mm]
{\it Philippe Bolle,
D\'epartement de math\'ematiques, Universit\'e
d'Avignon, 33, rue Louis Pasteur, 84000 Avignon, France,
philippe.bolle@univ-avignon.fr}

\end{document}

%% file: macro.tex
\newtheorem{theorem}{Theorem}[section]

\newtheorem{corollary}{Corollary}[section]
\newtheorem{condition}{Condition}[section]
\newtheorem{remark}{Remark}[section]
\newtheorem{remarks}{Remark}[section]
\newtheorem{definition}{Definition}[section]
\newcommand{\be}{\begin{equation}}
\newcommand{\ee}{\end{equation}}

\newcommand{\teta}{\theta}

\newcommand{\om}{\omega}

\newcommand{\ep}{\epsilon}

\newcommand{\ov}{\overline}

\newcommand{\wtilde}{\widetilde}

\renewcommand{\a }{\alpha }

\renewcommand{\d }{\delta }

\newcommand{\e }{\varepsilon }

\newcommand{\vphi}{\varphi }

\renewcommand{\t }{\tau }



%% file: lincei11900.bbl
\begin{thebibliography}{10}

\bibitem{AB}
A. Ambrosetti, M. Badiale,
{\it Homoclinics: Poincar\'e-Melnikov type results
via a variational approach},
C. R. Acad. Sci. Paris, t. 323, S\'erie I, 1996, 753-758, and
Annales I. H. P. - Analyse nonlin., vol. 15, n.2, 1998, p. 233-252.

\bibitem{An1} S. Angenent: {\it A variational interpretation
of  Melnikov's function and exponentially small separatrix splitting},
Lecture notes of the London Math. Soc, Symplectic geometry,
ed. Dietmar Salamon.

\bibitem{Arn} V. I. Arnold: {\it Instability of dynamical systems
with several degrees of freedom}, Sov. Math. Dokl. 6, 1964, p. 581-585.

\bibitem{BB} M. Berti, P. Bolle: {\it Homoclinics and Chaotic Behaviour
for Perturbed Second order Systems}, Annali di Mat. Pura e Applicata,
(IV), vol. CLXXVI, 1999, pp. 323-378.

\bibitem{BBS} M. Berti, P. Bolle: {\it Variational construction of
Homoclinics and Chaotic Behaviour in presence of a saddle-saddle
equilibrium}, Annali della Scuola Normale Superiore di Pisa,
serie IV, vol. XXVII, fasc. 2, 1998 and
Rend. Mat. Acc. Naz. Lincei, s. 9, vol. 9, fasc. 3, 1998.

\bibitem{BB1} M. Berti, P. Bolle: {\it Diffusion time and splitting
of separatrices for nearly integrable isochronous Hamiltonian systems},
to appear.

\bibitem{BCV} U. Bessi, L. Chierchia, E. Valdinoci:
{\it Lower Bounds on Arnold Diffusion Time via Mather theory}, preprint.


\bibitem{CG} L. Chierchia, G. Gallavotti: {\it Drift and diffusion
in phase space}, Annales de l'IHP, section Physique Th\'eorique, 60,
pp. 1-144, 1994; see also Erratum in Vol. 68, 135, 1998.


\bibitem{C1} J. Cresson: {\it Conjecture de Chirikov et
Optimalit\'e des exposants de stabilit\'e du th\'eor\`em
de Nekhoroshev}, preprint univ. Besancon.

\bibitem{DGTJ}  A. Delshams, V. G. Gelfreich,
V. G. Jorba, T. M. Seara: {\it Exponentially small
splitting  of separatrices under fast quasi-periodic forcing}
Comm. Math Ph. 189, 35-71, 1997.

\bibitem{G1}  G. Gallavotti: {\it Arnold's Diffusion in Isochronous
Systems}, Mathematical Physics, Analysis and Geometry 1, 295-312, 1999.



\bibitem{GGM}  G. Gallavotti, G. Gentile, V. Mastropietro:
{\it Separatrix splitting for systems with three time scales},
Commun. Math. Phys. 202, 197-236, 1999.


\bibitem{GGM3}  G. Gallavotti, G. Gentile, V. Mastropietro:
{\it A possible counter example to a paper by Rudnew and Wiggins},
Physica D, 137, 202-204, 2000.



\bibitem{Lo} P. Lochak
{\it Arnold diffusion: a compendium of remarks and questions},
Proceedings of 3DHAM's Agaro, 1995.


\bibitem{LMS} P. Lochak, J.P Marco, D. Sauzin,
{\it On the splitting of invariant
manifolds in multidimensional Hamiltonian systems}, preprint.

\bibitem{M} J. P. Marco {\it
Transitions le long des cha\^{\i}nes de tores invariants
pour les syst\`emes hamiltoniens analytiques},
Annales I. H. P., vol. 64, 1995, p. 205-252.


\bibitem{PV} A. Pumarino, C. Valls: {\it Three time scales systems
exhibiting persisent Arnold'Diffusion}, preprint.

\end{thebibliography}
